\documentclass[12pt]{amsart}
\usepackage[margin=0.5in]{geometry}
\newcommand{\R}{\mathbb{R}}
\newcommand{\Z}{\mathbb{Z}}
\newcommand{\N}{\mathbb{N}}
\newcommand{\Q}{\mathbb{Q}}

\renewcommand{\epsilon}{\varepsilon}
\newcommand{\eps}{\varepsilon}

\usepackage[english]{babel}
\usepackage[alphabetic]{amsrefs}
\usepackage{amsmath,amssymb,amsfonts,amsthm,enumerate}
\usepackage{url}
\usepackage{graphicx,epstopdf,color}

\numberwithin{equation}{section}

\newtheorem{theorem}{Theorem}[section]

\newtheorem{conjecture}[theorem]{Conjecture}

\renewcommand{\le}{\leqslant}

\renewcommand{\ge}{\geqslant}

\newenvironment{dedication}
        {\vspace{6ex}\begin{quotation}\begin{center}\begin{em}}
        {\par\end{em}\end{center}\end{quotation}\vspace{6ex}}

\begin{document}

\title[Nonlocal phase transitions]{Nonlocal phase transitions\\
in homogeneous and periodic media}

\author[Matteo Cozzi, Serena Dipierro, Enrico Valdinoci]{
Matteo Cozzi${}^{(1)}$
\and
Serena Dipierro${}^{(1,2)}$
\and
Enrico Valdinoci${}^{(1,2,3)}$
}

\maketitle

{\scriptsize \begin{center}(1) -- Weierstra{\ss}
Institut f\"ur Angewandte Analysis und Stochastik\\
Mohrenstra{\ss}e 39, D-10117 Berlin (Germany).
\end{center}
\scriptsize \begin{center} (2) -- 
School of Mathematics and Statistics  \\
University of Melbourne\\
Grattan Street, 
Parkville, VIC-3010 Melbourne
(Australia).\\
\end{center}
\scriptsize \begin{center} (3) -- Dipartimento di Matematica
``Federigo Enriques''\\
Universit\`a
degli studi di Milano\\
Via Saldini 50, I-20133 Milano (Italy).\\
\end{center}
\bigskip

\begin{center}
E-mail addresses: matteo.cozzi@wias-berlin.de,
serydipierro@yahoo.it,
enrico@math.utexas.edu
\end{center}
}

\begin{dedication}
\hspace{1cm}
{In honor of Professor 
Paul Rabinowitz, with great esteem and admiration.}
\end{dedication}

\begin{abstract}
We discuss some results related to a phase transition model
in which the potential energy induced by a double-well function
is balanced by a fractional elastic energy.
In particular, we present asymptotic results
(such as $\Gamma$-convergence, energy bounds and density estimates
for level sets), flatness
and rigidity results, and the construction of planelike
minimizers in periodic media.

Finally, we consider a nonlocal equation with a multiwell
potential, motivated by models arising in crystal dislocations,
and we construct orbits exhibiting symbolic dynamics,
inspired by some classical results by
Paul Rabinowitz.
\end{abstract}

\section{Introduction}

Goal of this paper is to collect into a homogeneous
and original presentation a series of results
about the fractional Allen-Cahn (or scalar
Ginzburg-Landau) equation.

The classical model for this equation arises
in the study of phase coexistence and it has several
applications in material sciences. In its basic
version, the model aims to describe the phase separation
occurring in some media. The two phases can be described
by a state parameter function~$u:\R^n \to [-1,1]$.
In this setting,~$n\in\N$ is the dimension
of the ambient space and the values~$-1$ and~$1$ for~$u$
represent the ``pure phases'' of the system.
\medskip

The total energy of the system is supposed to be made of
two terms: a ``potential'' term~${\mathcal{W}}(u)$,
which forces minimizers to stay ``as close as possible
to the pure phases'', and an ``elastic'' (or, with
a slight abuse of a terminology borrowed from
similal setting in Hamiltonian dynamics, ``kinetic'') term~${\mathcal{K}}(u)$,
which ``prevents the formation of unnecessary
interfaces''.\medskip

More precisely, given a (bounded, smooth) set~$\Omega\subset\R^n$
(the ``container''),
and a smooth function~$W:\R\to[0,+\infty]$
(the ``potential'')
which has nondegenerate minima at~$\pm1$ (i.e., say, $W(\pm1)=0<W(r)$
for any~$r\in\R\setminus(-1,1)$, with~$W''(\pm1)>0$), we set
$$ {\mathcal{W}}(u):=
\int_\Omega W(u(x))\,dx.$$
In the classical case, the formation of interfaces
is penalized by a local elastic term of the form
\begin{equation}\label{K1}
{\mathcal{K}}(u):=\frac12\int_\Omega |\nabla u(x)|^2\,dx,
\end{equation}
which leads to the total energy functional
\begin{equation}\label{E1}
{\mathcal{E}}(u):= {\mathcal{K}}(u)+
{\mathcal{W}}(u)=\int_\Omega \frac{|\nabla u(x)|^2}{2}+W(u(x))\,dx,
\end{equation}
whose critical points are solutions of the partial
differential equation
\begin{equation}\label{D1}
-\Delta u(x)+W'(u(x))=0, \qquad{\mbox{ for any }}x\in\Omega.
\end{equation}

Recently, some attention has been devoted
to nonlocal phase transition models, in order to
capture the long term interactions between particles
and to describe the boundary effects, see e.g.~\cite{BELL, MAR, SIRE_VAL_IFB}.
Here, we consider a phase transition model driven
by a nonlocal energy of fractional type, which can be described
as follows. Particles are supposed to interact
according to a kernel, which we take invariant
under translations and rotations, scale invariant
and with polynomial decay. More concretely, we set
\begin{equation} \label{Kdef}
K(y):=\frac{1}{|y|^{n+2s}},
\end{equation}
with~$s\in(0,1)$ and consider, as elastic energy, the quantity
\begin{equation}\label{Ks}
{\mathcal{K}}(u):=\frac12\int_{Q_\Omega} |u(x)-u(y)|^2\;K(x-y)\,dx\,dy,\end{equation}
where
$$ Q_\Omega:= \R^{2n}\setminus (\Omega^c)^2
=(\Omega\times\Omega)\cup (\Omega\times \Omega^c)
\cup (\Omega^c\times \Omega).$$
Here above, $\Omega^c:=\R^n\setminus\Omega$ denotes the complementary set.

By comparing~\eqref{K1} with~\eqref{Ks},
we see that we have replaced the classical seminorm in the Sobolev
space~$H^1$ with a seminorm in the fractional Sobolev space~$H^s$,
with~$s\in(0,1)$. As for the domain of integration,
the idea, both in~\eqref{K1} and in~\eqref{Ks},
is that the values of the state parameter~$u$
are fixed outside the container~$\Omega$,
so they should not really contribute to an effective energy
and the energy should only take into account contributions
which ``see the
container~$\Omega$''.
In this sense, the integral in~\eqref{K1}
takes into account the whole of the space~$\R^n$,
with the exception of the contributions
that lie outside the container~$\Omega$
(that is, the integral in~\eqref{K1} ranges in~$\R^n \setminus(\Omega^c)=
\Omega$). In the same spirit, the energy in~\eqref{Ks},
which is a double integral, takes into account all
the interactions in the whole of the space~$\R^n\times\R^n$,
with the exception of the ones
which only involve points outside the container~$\Omega$
(that is, the integral in~\eqref{Ks} ranges in~$(\R^n\times\R^n)
\setminus(\Omega^c\times\Omega^c)$, which indeed coincides with~$Q_\Omega$).
\medskip

With these observations, when the elastic energy in~\eqref{K1}
is replaced by the nonlocal one in~\eqref{Ks}, the total
energy in~\eqref{E1} is replaced by its fractional analogue
\begin{equation}\label{Es}
{\mathcal{E}}(u):= {\mathcal{K}}(u)+
{\mathcal{W}}(u)=
\frac12\int_{Q_\Omega} |u(x)-u(y)|^2\;K(x-y)\,dx\,dy
+\int_\Omega W(u(x))\,dx,
\end{equation}
whose critical points are solutions of the partial
differential equation
\begin{equation}\label{Ds}
(-\Delta)^s u(x)+W'(u(x))=0, \qquad{\mbox{ for any }}x\in\Omega,
\end{equation}
which can be seen as a nonlocal counterpart of~\eqref{D1}.
Here, as customary in the literature involving
fractional operators, we are using the notation~$(-\Delta)^s$
to denote the fractional Laplacian, i.e.~the integrodifferential
operator given (up to multiplicative dimensional constants
that we neglect here) by
$$ \int_{\R^n} \big( 2u(x)-u(x+z)-u(x-z)\big)\,K(z)\,dz.$$
We refer to~\cite{Land, Stein, Silv-TH, guida} for an introduction
to the fractional Laplacian.\medskip

Here, we aim to discuss the theory
of nonlocal phase transitions, as described by
the energy functional in~\eqref{Es} and by the
pseudodifferential equation in~\eqref{Ds}, basically
discussing the similarities
and the important differences with respect to the classical
theory, especially in the light of
$\Gamma$-convergence, density estimates,
rigidity and flatness results and periodic and
quasiperiodic structures
arising in periodic media.

\section{$\Gamma$-convergence and density estimates}

To study the asymptotics of the solutions
of~\eqref{Ds}, it is convenient to look
at the spacial scaling~$x\to\frac{x}{\eps}$.
Correspondingly, one can appropriately
rescale the energy functional as
$$ {\mathcal{E}}_\eps(u):=
\begin{cases}
{\mathcal{K}}(u) + \eps^{-2s} {\mathcal{W}}(u)
& {\mbox{ if }} s\in (0,\,1/2),\\
|\log\eps|^{-1}
{\mathcal{K}}(u) + |\eps\log\eps|^{-1} {\mathcal{W}}(u)
& {\mbox{ if }} s=1/2,\\
\eps^{2s-1}
{\mathcal{K}}(u) + \eps^{-1} {\mathcal{W}}(u)
& {\mbox{ if }} s\in(1/2,\,1).
\end{cases} $$
When we want to emphasize the dependence of the energy
functional (or of the energy contributions)
on the container~$\Omega$,
we will write, respectively, ${\mathcal{W}}(u;\Omega)$, 
${\mathcal{K}}(u;\Omega)$, ${\mathcal{E}}(u;\Omega)$
and~${\mathcal{E}}_\eps(u;\Omega)$.

In this notation, if~$u_\eps(x):=u(x/\eps)$ and~$\Omega_\eps:=
\eps\Omega$, we have that
\begin{eqnarray*}&&
{\mathcal{W}}(u_\eps;\Omega_\eps) = \eps^n \,{\mathcal{W}}(u;\Omega)
\\ {\mbox{and }}&&
{\mathcal{K}}(u_\eps;\Omega_\eps) = \eps^{n-2s} \,{\mathcal{K}}(u;\Omega)
\end{eqnarray*}
and so
$$ {\mathcal{E}}_\eps(u_\eps;\Omega_\eps)=
\eps^{n-\min\{2s,1\}} {\mathcal{E}}(u;\Omega)$$
if~$s\in(0,1)\setminus\{1/2\}$, and
$$ {\mathcal{E}}_\eps(u_\eps;\Omega_\eps)=\frac{
\eps^{n-1} }{|\log\eps|}{\mathcal{E}}(u;\Omega)$$
if~$s=1/2$ (this additional logaritmic factor is indeed
very typical for fractional problems related to the square root
of the Laplacian).\medskip

The scale of this functional is chosen in such a way
that the following $\Gamma$-convergence result holds
(see~\cite{SV-Gamma}):

\begin{theorem} \label{GAMMA:CO}
As~$\eps\searrow0$, the functional ${\mathcal{E}}_\eps$
$\Gamma$-converges to
$$ {\mathcal{E}}_0(u):=
\begin{cases}
{\mathcal{K}}(u;\Omega) 
& {\mbox{ if }} s\in (0,\,1/2) {\mbox{ and }}
u=\chi_E-\chi_{E^c} {\mbox{ a.e. in }} \Omega {\mbox{ for some }}
E\subseteq\Omega,\\
c_\star\, {\rm Per}\,(E,\Omega)
& {\mbox{ if }} s\in [1/2,\,1) {\mbox{ and }}
u=\chi_E-\chi_{E^c} {\mbox{ a.e. in }} \Omega {\mbox{ for some }}
E\subseteq\Omega,\\
+\infty &{\mbox{ otherwise.}}
\end{cases} $$
\end{theorem}

It is worth to point out that Theorem~\ref{GAMMA:CO}
may be rephrased by saying that
when~$s\in [1/2,\,1)$ the fractional phase transition energy
$\Gamma$-converges to the classical perimeter functional
-- as it happens indeed in the classical case~$s=1$.
The classical counterpart when~$s=1$ of Theorem~\ref{GAMMA:CO}
was indeed one of the founding results of $\Gamma$-convergence,
see~\cite{mortola}.\medskip

On the other hand,
when~$s\in (0,\,1/2)$, the fractional phase transition energy
converges to the fractional
perimeter functional that was introduced in~\cite{ROQ}.
This suggests that the fractional parameter~$s=1/2$
provides a threshold for the asymptotic behavior of
nonlocal phase transitions: above such threshold,
the behavior of the interfaces at a large scale somehow
``localizes'', since such behavior is dictated by
the classical, and thus local, perimeter 
minimization; but below such threshold
the behavior of the interfaces at a large scale fully
mantains its nonlocal properties, since
it is driven by a nonlocal perimeter functional.\medskip

Of course, $\Gamma$-convergence is a very elegant and effective
method to deal with the asymptotics of functionals and it
fits well with the calculus of variations and with the
problems of minimization. On the other hand, it gives
little information on the geometric properties of 
the solutions of the equation. In the classical case, to
overcome this difficulty, a theory of ``density estimates''
has been developed in~\cite{cordoba}. Namely,
the goal of this theory is to establish energy bounds
and bounds in measure theoretic sense for the level sets of
minimizers, with the goal of showing that minimizers behave
``like one-dimensional solutions'' at least in terms of energy
and in terms of measure occupied by their interface.
\medskip

In the fractional framework, we have that
the energy of minimizers of~${\mathcal{E}}_\eps$
is locally uniformly bounded, according to the following result:

\begin{theorem}\label{EBO}
Let~$\vartheta_1\in(0,1)$.
If $u_\eps$ minimizes~${\mathcal{E}}_\eps$
in $B_{1+\eps}$ and~$|u_\eps(0)|<\vartheta_1$, then
\begin{equation}\label{EN:EST}
c\le {\mathcal{E}}_\eps (u_\eps,B_1)\le C,\end{equation}
with~$C>c>0$ only depending on $n$, $s$ and~$ W$.
\end{theorem}

The upper energy bound in~\eqref{EN:EST} was proved in~\cite{SV-DENS} 
and the lower bound in~\cite{cozzi-valdinociB}.\medskip 

A counterpart of the energy bounds in Theorem~\ref{EBO}
is a density estimate,
which says, roughly speaking, that the measure of
the interface of minimizers of~${\mathcal{E}}_\eps$
is locally of size comparable to~$\eps$. The precise result
goes as follows:

\begin{theorem}\label{DENS}
Let~$\vartheta_1$, $\vartheta_2\in(0,1)$.
If $u_\eps$ minimizes~${\mathcal{E}}_\eps$
in $B_{1}$ and~$|u_\eps(0)|<\vartheta_1$, then
\begin{equation}\label{DENS:EST}
c\eps\le \big|\{ |u_\eps|<\vartheta_2 \}\cap B_1 \big|\le C\eps,\end{equation}
with~$C>c>0$ only depending on $n$, $s$ and~$ W$.
\end{theorem}

The upper density estimate in~\eqref{DENS:EST} was proved in~\cite{SV-DENS}
and the lower bound in~\cite{cozzi-valdinociB}.
We remark that the connection between the energy bound
in~\eqref{EN:EST} and the measure theoretic bound in~\eqref{DENS:EST}
is particularly close in the local case (i.e.~for~$s=1$)
and also in the weakly nonlocal case (i.e. when~$s\in(1/2,\,1)$),
since, roughly speaking, the potential energy is capable
of measuring the size of the interface in a sufficiently
sharp way. On the other hand, in
the strongly 
nonlocal case (i.e.~when~$s\in(0,\,1/2)$),
the potential energy does not suffice for this scope
and one has to carefully take into account the interaction
in the elastic energy, possibly at any scale, to detect
the dominant contributions. In particular, the case~$s=1/2$
for these estimates turns out to be the most delicate
one, since the precise bounds involve a logaritmic correction.
\medskip

We stress that the optimal bounds
obtained in~\eqref{EN:EST} and~\eqref{DENS:EST}
not only provide the uniform convergence of the level sets
of the minimizers to their limit interface (see~\cite{cordoba}),
but also provide the cornerstone for the construction
of planelike solutions in periodic media for which the oscillation
from the reference plane is of the same order of the size of
periodicity of the media (this feature will be more
throughly detailed in Section~\ref{PLANE:SEC}). 
For this, we also remark that both Theorems~\ref{EBO} 
and~\ref{DENS} hold for more general kernels~$K$ and 
potentials~$W$ (see~\cite{cozzi-valdinociB}).

\section{Rigidity and flatness results} \label{RIGFLAT:SEC}

A byproduct of the optimal bounds
in~\eqref{EN:EST} and~\eqref{DENS:EST} is that minimizers
behave as if they were one-dimensional functions
(that is, for functions of only one Euclidean variable,
with sufficient decay at infinity, one can check ``by hands'' that
formulas~\eqref{EN:EST} and~\eqref{DENS:EST} hold true).
\medskip

A natural question in this setting is whether global
one dimensional solutions of~\eqref{Ds} indeed exist. That is,
up to normalization, if there exists a function~$u_0:\R\to(-1,1)$
that satisfies
\begin{equation}
\label{trale}\begin{split} &
(-\Delta)^s u_0(x)+W'(u_0(x))=0, \qquad{\mbox{ for any }}x\in\R,
\\&\lim_{x\to\pm\infty} u_0(x) =\pm1.
\end{split}\end{equation}
We stress that when~$s=1$, the existence of such ``transition layer''~$u_0$
is obvious, since the problem boils down to an ordinary
differential equation, which can be integrated explicitly
(by multiplying the equation by~$u_0'(x)$ and taking
the antiderivative -- or equivalently by using the
Law of Conservation of Energy).\medskip

On the other hand, differently from the classical case,
when~$s\in(0,1)$ the existence of such solution~$u_0$
is a rather delicate business, and it has been
established, using variational and energy methods, in~\cite{PSV-AMPA},
\cite{CS-Trans} and, in further generality, in~\cite{tommi}.\medskip

Of course, given the estimates in~\eqref{EN:EST} and~\eqref{DENS:EST},
one may wonder under which additional conditions (if any)
we can say that solutions of~\eqref{Ds}
are indeed one-dimensional, i.e., up to normalizations,
are of the form~$u(x)=u_0(x_n)$, or, say,
the level sets of~$u$ are hyperplanes. In the classical case~$s=1$
this was in fact the content of a beautiful conjecture by Ennio De Giorgi
in~\cite{DG-CONJ}, which can be stated as follows:

\begin{conjecture}\label{CON-DG79}
Let~$u\in C^2(\R^n,[−1,1])$
satisfy
\begin{equation}\label{DG:EQ}
-\Delta u= u-u^3
\end{equation}
and
$$ \partial_{x_n} u>0$$
in the whole of~$\R^n$.
Is it true that all the level sets of~$u$
are hyperplanes, at least if~$n\le8$?
\end{conjecture}

We refer to~\cite{STATEOFTHEART} for a detailed
account of the available results related to
Conjecture~\ref{CON-DG79};
here, we would like to discuss the fractional
analogue of Conjecture~\ref{CON-DG79} when equation~\eqref{DG:EQ}
is replaced by its nonlocal counterpart
\begin{equation*}
(-\Delta)^s u= u-u^3
\end{equation*}
with~$s\in(0,1)$. In this case,
a positive answer to this problem was given in~\cite{CSM-CPAM}
when~$n=2$ and~$s=1/2$, in~\cite{SireV-JFA, CS-Trans}
when~$n=2$ and~$s\in(0,1)$, in~\cite{CC-DCDS}
when~$n=3$ and~$s=1/2$ and in~\cite{CC-CalcVar}
when~$n=3$ and~$s\in(1/2,\,1)$
(see also~\cite{SV-JFA} and~\cite{bucur-monogr}
for different proofs, also related to nonlocal minimal surfaces).
\medskip

Hence, all in all, at the moment, to the best of our knowledge,
the state of the art on the nonlocal analogue of
Conjecture~\ref{CON-DG79}
can be summarized by the following result:

\begin{theorem}\label{HGA:TH}
Let~$u\in C^2(\R^n,[−1,1])$
satisfy
\begin{eqnarray*}
&& (-\Delta)^s u= u-u^3
\\ {\mbox{and }}&&\partial_{x_n} u>0\end{eqnarray*}
in the whole of~$\R^n$.
Assume also that
\begin{equation}\label{COsATT}
\begin{split}
&{\mbox{ either $n=2$ and $s\in(0,1)$,}}\\
&{\mbox{ or $n=3$ and $s\in [1/2,\,1)$.}}
\end{split}
\end{equation}
Then all the level sets of~$u$
are hyperplanes.
\end{theorem}

As a matter of fact, Theorem~\ref{HGA:TH} holds true
for a very general class of equations under condition~\eqref{COsATT},
see~\cite{CSM-CPAM, SireV-JFA, CS-Trans, CC-DCDS, CC-CalcVar, SV-JFA, bucur-monogr}.
Of course, it would be very desirable to go beyond
condition~\eqref{COsATT}, or to provide counterexamples.\medskip

It is also suggestive to compare the threshold~$s=1/2$
in~\eqref{COsATT} with the different behavior of the $\Gamma$-limit
described by Theorem~\ref{GAMMA:CO}. In any case,
it is not clear whether or not condition~\eqref{COsATT}
reflects somehow the different behavior of local and nonlocal
minimal surfaces. See also~\cite{bucur-monogr}
for further discussions on this point.

\section{Planelike minimizers in periodic media}\label{PLANE:SEC}

A classical topic in differential geometry (resp.,
dynamical systems) is to look for periodic and
quasiperiodic solutions of problems set in
periodic media which lie at a bounded distance
from any fixed hyperplane (resp. which possess a rotation number
in average). For instance, in~\cite{MORSE, HEDL}
the case of the plane with a Riemannian metric was
taken into account, establishing the existence
of minimal geodesics which stay at a bounded
distance from any prescribed straight line
(depending on the arithmetic 
properties of the slope of this line,
the geodesics turn out to be either periodic
or quasiperiodic).\medskip

A similar problem in higher codimension turns out
to be, in general, ill-posed, since~\cite{HEDL}
provided an example of a metric in~$\R^3$
which does not have geodesics at bounded distance
from a particular direction. Nevertheless,
a similar problem can be efficiently set
in higher dimension, provided that one 
takes into consideration objects of codimension~$1$, such as
minimal hypersurfaces (namely, hypersurfaces
which minimize a periodic perimeter functional)
that lie at bounded distance from any fixed
hyperplane. We refer to~\cite{Bangert90, CLL01, Auer01}
and references therein for these types of problems
and, for instance, to~\cite{Mather}
for related problems in the setting of dynamical systems.
\medskip

The construction of periodic and quasiperiodic objects
in periodic media has also a long and important
tradition in partial differential equations, see~\cite{Moser86, RaStre-book}
and the references therein.
For instance, in~\cite{V-CRELLE} the classical
Ginzburg-Landau-Allen-Cahn equation~\eqref{DG:EQ}
is set in a periodic medium, and one constructs
global solutions which are minimal, and either
periodic or quasiperiodic, and whose level
sets stay at a bounded distance from a prescribed
hyperplane. Roughly speaking, from a physical point of
view, one obtains in this way
some phase coexistence in an infinite periodic
medium whose interface is a flat hyperplane, up to a uniformly
bounded error. \medskip

Goal of this section is to describe in detail some
of the results concerning periodic and quasiperiodic
minimizers for the nonlocal phase transition
equation in a periodic medium. For this, we take into consideration the following heterogeneous variant of the total energy~\eqref{Es}:
\begin{equation} \label{Epl}
{\mathcal{E}}(u) = {\mathcal{E}}(u; \Omega) := \frac{1}{2} \int_{Q_\Omega} |u(x) - u(y)|^2 \; K(x - y) \, dx \, dy + \int_{\Omega} Q(x) \, W(u(x)) \, dx,
\end{equation}
with~$K$ as in~\eqref{Kdef}, and the corresponding Euler-Lagrange equation
$$
(-\Delta)^s u(x) + Q(x) \, W'(u(x)) = 0, \qquad \mbox{ for any } x \in \R^n.
$$

In order to model a non-homogeneous periodic environment, 
the potential~$W$ appears now modulated by a measurable 
function~$Q: \R^n \to [Q_*, Q^*]$, with~$Q^* \ge Q_* > 0$, 
which is periodic with respect to a discrete lattice of step~$\tau \ge 1$, that is
$$
Q(x + k) = Q(x), \qquad \mbox{ for a.e. } x \in \R^n \mbox{ and any } k \in \tau \Z^n.
$$

\smallskip

As we just said, the energy~$\mathcal{E}$ embodies the presence of an underlying heterogeneous medium by means of the multiplicative correction~$Q$ in the potential term. Note that the model is sensitive to the periodicity scale of the medium through the factor~$\tau$.

Of course, a broader setting can be taken into account, for instance by considering a more general periodic potential or by
letting the kernel~$K$ be space-dependent as well. 
This latter generalization is particularly interesting as it encompasses models in which the interaction between two particles of the system does not depend only on their distance, but may vary (albeit in a periodic way) as the particles occupy different places in the space.

Here, we choose to favor a not too involved exposition and thus to stick to the simpler model provided by~\eqref{Epl}. For a presentation of our contributions in a wider generality, we refer the interested reader to~\cite{cozzi-valdinociA, cozzi-valdinociB}.
\medskip

The result that we shortly discuss addresses the existence of a~\emph{planelike} minimizer for~$\mathcal{E}$ in the whole space~$\R^n$. That is, a function~$u: \R^n \to [-1, 1]$ that minimizes~$\mathcal{E}$ in~$\Omega$ for any bounded set~$\Omega \subset \R^n$ (a so-called~\emph{class~A minimizer}) and whose intermediate level set~$\{ |u| < 9 / 10 \}$ lies at a bounded distance from any fixed hyperplane of~$\R^n$. By construction, the minimizer enjoys a suitable periodicity or quasi-periodicity property, depending on whether the slope of the associated hyperplane is rational or not.

When~$Q$ is constant, the energy is translation-invariant and, as a result,~\emph{planar} minimizers exist (recall the discussion at the beginning of Section~\ref{RIGFLAT:SEC}). Under the presence of a non-trivial modulation~$Q$, such construction cannot be carried out and, in general, one-dimensional minimizers do not exist. However, the following fact holds true.

\begin{theorem} \label{PLthm}
There exists a universal constant~$M_0 > 0$, such that, given any~$\omega \in \R^n \setminus \{ 0 \}$, we can construct a class~A minimizer~$u$ for~$\mathcal{E}$ satisfying
$$
\left\{ x \in \R^n : |u(x)| < \frac{9}{10} \right\} \subset \left\{ x \in \R^n : \frac{\omega}{|\omega|} \cdot x \in [0, \tau M_0] \right\}.
$$
Furthermore,~$u$ enjoys the following quasi-periodicity properties:
\begin{enumerate}[$\, \bullet$]
\item if~$\omega \in \tau \Q^n \setminus \{ 0 \}$, then~$u$ is~$\sim$-periodic, 
i.e.~it respects the equivalence relation~$\sim$ defined in~$\R^n$ by
\begin{equation} \label{tildef}
x \sim y \qquad \mbox{ if and only if } \qquad y - x \in \tau \Z^n \mbox{ and } \omega \cdot (y - x) = 0;
\end{equation}
\item if~$\omega \in \R^n \setminus \tau \Q^n$, then~$u$ is the locally uniform limit of a sequence of periodic class~A minimizers.
\end{enumerate}
\end{theorem}

\begin{figure}
    \centering
    \includegraphics[width=16.8cm]{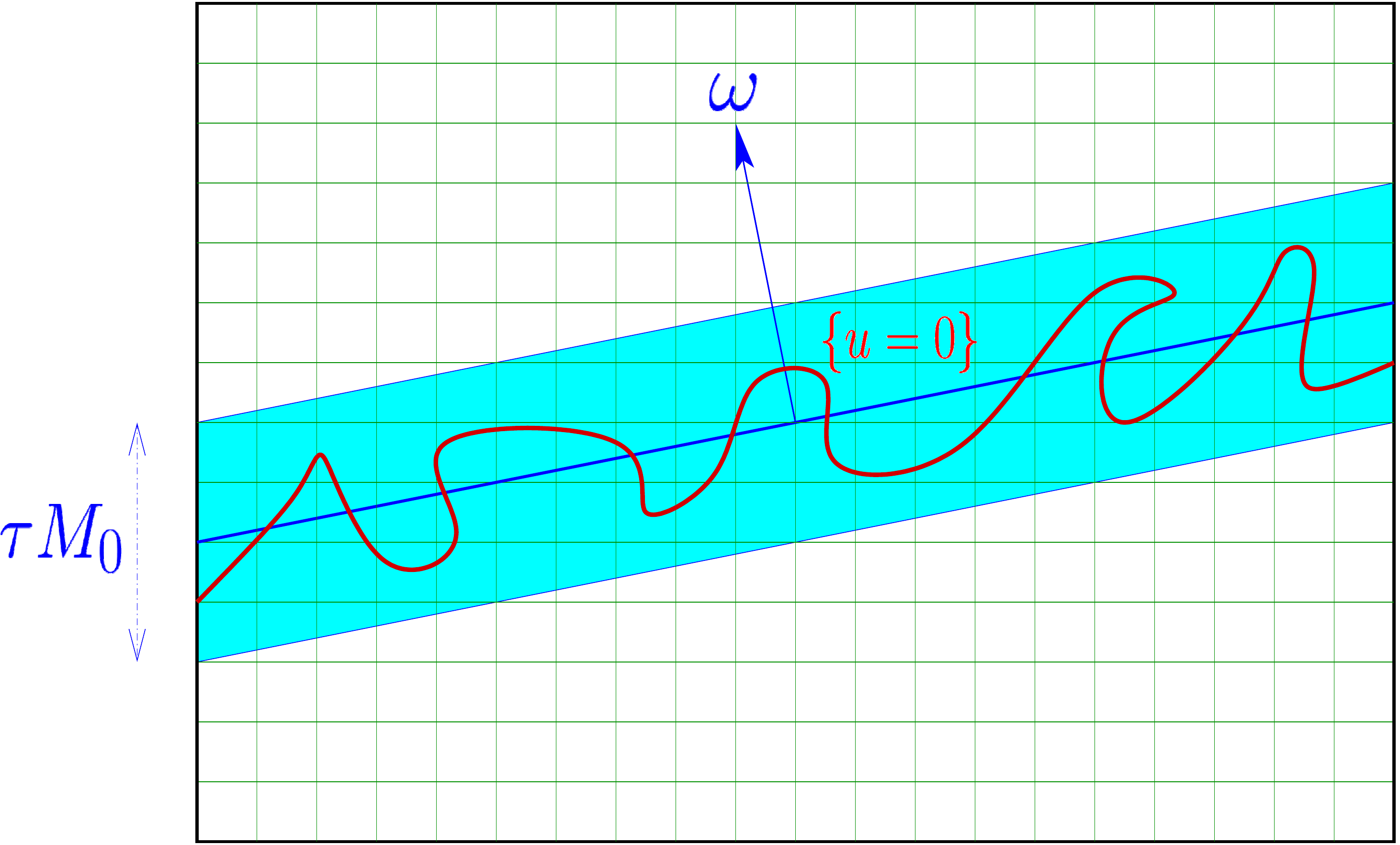}
    \caption{A planelike minimizer as given in Theorem~\ref{PLthm}.}
    \label{fig:plane}
\end{figure}

Theorem~\ref{PLthm} provides the aforementioned 
existence of planelike minimizers for the nonlocal energy~\eqref{Epl}, 
see Figure~\ref{fig:plane}.

We stress that the size of the strip where the (essential) transition of a minimizer occurs is proportional to the periodicity scale~$\tau$ of the medium. The proportionality constant~$M_0$ is universal, in the sense that it depends only on the structural constants involved in the model, i.e.~$n$,~$s$,~$Q_*$ and~$Q^*$.

Besides being interesting in itself, this fact plays a crucial role in deducing from Theorem~\ref{PLthm} a similar construction of planelike minimal surfaces for a periodic nonlocal perimeter (see~\cite{cozzi-valdinociB} for more details).

Also, the value~$9/10$, that has been used to identify an interface region for the minimizer, obviously plays no particular role and may be indeed replaced by any~$\vartheta \in (0, 1)$. However, the constant~$M_0$ would then depend on~$\vartheta$ as well.
\medskip

The rest of this section contains a sketch of the proof of Theorem~\ref{PLthm}. 
The general strategy adopted is shaped on the 
one designed in~\cite{CLL01} for a model described 
by a periodic surface energy, and developed in~\cite{V-CRELLE}, 
where the argument is used to deal with a local, non-homogeneous 
Allen-Cahn-Ginzburg-Landau functional similar to~\eqref{E1}.
\smallskip

We begin by addressing the case of a 
direction~$\omega \in \tau \Q^n \setminus \{ 0 \}$. 
We denote by~$\widetilde{\R}^n$ any fundamental domain of 
the quotient space~$\R^n / \sim$, where the relation~$\sim$ is given in~\eqref{tildef}, 
and consider the class of admissible functions
$$
\mathcal{A}_\omega^M := \left\{ u \in L^2_{\rm loc}(\R^n) : u \mbox{ is~$\sim$-periodic, } u(x) \ge \frac{9}{10} \mbox{ if } \frac{\omega}{|\omega|} \cdot x \le 0 \mbox{ and } u(x) \le - \frac{9}{10} \mbox{ if } \frac{\omega}{|\omega|} \cdot x \ge M \right\},
$$
for any large~$M > 0$. A straightforward application of the 
Direct Method of the Calculus of Variations gives the existence of global minimizers of the auxiliary functional
$$
{\mathcal{F}}(u) := \frac{1}{2} \int_{\widetilde{\R}^n} 
\int_{\R^n} |u(x) - u(y)|^2 \; K(x - y) \, dx \, dy + \int_{\widetilde{\R}^n} Q(x) \, W(u(x)) \, dx,
$$
within the set~${\mathcal{A}}_\omega^M$, 
at least\footnote{We stress that the case in which~$s \in (0, 1/2]$ 
can be treated afterwards via a limiting argument, approximating~$K$ 
with kernels truncated at infinity. The difficulties arise essentially 
for the fact that, when~$s \in (0,1/2]$, the 
functional~$\mathcal{F}$ is identically equal to~$+\infty$ 
on~${\mathcal{A}}_\omega^M$, due to the~``fat'' tails of the kernel.} 
if~$s \in (1/2, 1)$. So we denote by~${\mathcal{M}}_\omega^M$ the set of such minimizers.

Observe that~$\mathcal{F}$ differs from the 
energy~${\mathcal{E}}(\cdot, \widetilde{\R}^n)$, given in~\eqref{Epl}, 
for the sole fact that the latter contains the term integrated 
over~$\widetilde{\R}^n \times {( \widetilde{\R}^n )}^c$ counted twice, namely
\begin{equation}\label{UA:89A0101AL}
\begin{split}
Q_{\widetilde{\R}^n }& \;=\; \big( {\widetilde{\R}^n }\times{\widetilde{\R}^n }\big)\cup
\big( {\widetilde{\R}^n }\times({\widetilde{\R}^n })^c\big)\cup
\big( ({\widetilde{\R}^n })^c\times{\widetilde{\R}^n }\big)\\
{\mbox{while }}\qquad
{\widetilde{\R}^n } \times\R^n &\;=\;
\big( {\widetilde{\R}^n }\times{\widetilde{\R}^n }\big)\cup 
\big( {\widetilde{\R}^n }\times({\widetilde{\R}^n })^c\big).
\end{split}\end{equation}
The necessity of considering the auxiliary functional~$\mathcal{F}$ is 
peculiar to the nonlocal setting considered here and is mostly due to 
the fact that the functional~$\mathcal{E}$ does not behave 
well with respect to the periodic structure induced by~$\sim$. 

In more concrete terms, given a function~$u$, one can consider 
its~$\sim$-periodic extension~$\widetilde{u}$, defined for a.e.~$x \in \R^n$ as
$$
\widetilde{u}(x) := u(\widetilde{x}), \qquad \mbox{ where } 
\widetilde{x} \mbox{ is the only element of } \widetilde{\R}^n \mbox{ such that } 
x \sim \widetilde{x}.
$$
Then, it holds in general that
\begin{equation} \label{EneE}
{\mathcal{E}} (u; \widetilde{\R}^n) \ne {\mathcal{E}} (\widetilde{u}; \widetilde{\R}^n).
\end{equation}
However,
from \eqref{UA:89A0101AL}, noticing that the potential term has a local character,
we obtain that
\begin{equation}\label{PERO1}
\begin{split}
&{\mathcal{E}} (u;\widetilde{\R}^n) - {\mathcal{F}} (\widetilde u)
\\ =\;&
\frac12\int_{Q_{\widetilde{\R}^n}} |u(x)-u(y)|^2\;K(x-y)\,dx\,dy
-\frac{1}{2} \int_{\widetilde{\R}^n}
\int_{\R^n} |\widetilde u(x) - \widetilde u(y)|^2 \; K(x - y) \, dx \, dy\\ =\;&
\int_{({\widetilde{\R}^n })^c}
\int_{{\widetilde{\R}^n }} |u(x)-u(y)|^2\;K(x-y)\,dx\,dy
-\frac12
\int_{({\widetilde{\R}^n })^c}
\int_{{\widetilde{\R}^n }}
|\widetilde u(x) - \widetilde u(y)|^2 \; K(x - y) \, dx \, dy
.\end{split}
\end{equation}
In particular, if $u$ is $\sim$-periodic (hence $u=\widetilde u$),
\begin{equation}\label{PERO2}
{\mathcal{E}} (u;\widetilde{\R}^n) - {\mathcal{F}} (\widetilde u)=
\frac12\int_{({\widetilde{\R}^n })^c}
\int_{{\widetilde{\R}^n }} |u(x)-u(y)|^2\;K(x-y)\,dx\,dy
.\end{equation}
In addition, we have that
\begin{equation}\label{9iAJJ1A}\begin{split}
&{\mbox{any 
minimizer~$u$ of~$\mathcal{F}$ in the class~${\mathcal{A}}_\omega^M$ 
(i.e.~any~$u \in {\mathcal{M}}_\omega^M$)}}\\
&{\mbox{is a minimizer of the energy~$\mathcal{E}$ with respect to perturbations }}
\\&{\mbox{supported
inside the quotiented strip}}\\
&{\widetilde{\mathcal{S}}}_\omega^M := \widetilde{\R}^n \cap {\mathcal{S}}_\omega^M, \quad \mbox{ with } \quad {\mathcal{S}}_\omega^M := \left\{ x \in \R^n : \frac{\omega}{|\omega|} \cdot x \in [0, M] \right\}.
\end{split}\end{equation}
This fact, which is of key importance and ultimately motivates the introduction of the functional~$\mathcal{F}$,
is based on the following computation: let $\phi$ be a smooth function supported
inside ${\widetilde{\mathcal{S}}}_\omega^M$ and let $v:=u+\phi$, with $u\in{\mathcal{A}}_\omega^M$.
Then, $\phi$ vanishes outside ${\widetilde{\R}^n} $
and $\widetilde{v}=\widetilde{u}+\widetilde{\phi}=u+\widetilde{\phi}$. Consequently,
if $x\in ({\widetilde{\R}^n })^c$ and $y\in{\widetilde{\R}^n }$, we have that
\begin{eqnarray*}&&
|v(x)-v(y)|^2 - |\widetilde v(x)-\widetilde v(y)|^2\\
&=& |u(x)-u(y)-\phi(y)|^2 - |u(x)+\widetilde\phi(x)-u(y)-\phi(y)|^2\\
&=& |u(x)-u(y)|^2 +\phi^2(y)-2(u(x)-u(y))\phi(y)
\\ &&\qquad\quad -|u(x)+\widetilde\phi(x)-u(y)|^2-\phi^2(y)
+2(u(x)+\widetilde\phi(x)-u(y))\phi(y) \\
&=&|u(x)-u(y)|^2 
-|u(x)+\widetilde\phi(x)-u(y)|^2
+2\widetilde\phi(x)\,\phi(y).
\end{eqnarray*}
Hence, recalling \eqref{PERO1} and \eqref{PERO2},
\begin{equation}\label{TYUAG789Aaasfd}\begin{split}
& \Big({\mathcal{E}}(v;\widetilde{\R}^n)-{\mathcal{F}}(\widetilde v)\Big)
- \Big({\mathcal{E}}(u;\widetilde{\R}^n)-{\mathcal{F}}(u)\Big)
\\
=\;& 
\int_{({\widetilde{\R}^n })^c}
\int_{{\widetilde{\R}^n }} \left(|v(x)-v(y)|^2
-\frac12 |\widetilde v(x) - \widetilde v(y)|^2 
-\frac12 |u(x) - u(y)|^2 
\right)\; K(x - y) \, dx \, dy\\
=\;&
\int_{({\widetilde{\R}^n })^c}
\int_{{\widetilde{\R}^n }} \left(\frac12|v(x)-v(y)|^2
-\frac12|u(x)+\widetilde\phi(x)-u(y)|^2
+\widetilde\phi(x)\,\phi(y)
\right)\; K(x - y) \, dx \, dy
.\end{split}\end{equation}
Also, using the $\sim$-periodicity and changing variables $x\longmapsto x-k$,
$y\longmapsto y+k$,
\begin{eqnarray*}
&& \int_{({\widetilde{\R}^n })^c}
\int_{{\widetilde{\R}^n }}
|u(x)+\widetilde\phi(x)-u(y)|^2\; K(x - y) \, dx \, dy\\
&=& \sum_{{k\in\Z^n\setminus\{0\}}\atop{\omega\cdot k=0}}
\int_{{\widetilde{\R}^n }+k}
\int_{{\widetilde{\R}^n }}
|u(x)+\widetilde\phi(x)-u(y)|^2\; K(x - y) \, dx \, dy\\
&=&\sum_{{k\in\Z^n\setminus\{0\}}\atop{\omega\cdot k=0}}
\int_{{\widetilde{\R}^n }}
\int_{{\widetilde{\R}^n }-k}
|u(x)+\widetilde\phi(x)-u(y)|^2\; K(x - y) \, dx \, dy
\\
&=&
\int_{{\widetilde{\R}^n }}
\int_{({\widetilde{\R}^n })^c}
|u(x)+\widetilde\phi(x)-u(y)|^2\; K(x - y) \, dx \, dy
\\
&=&
\int_{{\widetilde{\R}^n }}
\int_{({\widetilde{\R}^n })^c}
|u(x)+\phi(x)-u(y)|^2\; K(x - y) \, dx \, dy
\\
&=&
\int_{{\widetilde{\R}^n }}
\int_{({\widetilde{\R}^n })^c}
|v(x)-v(y)|^2\; K(x - y) \, dx \, dy
.\end{eqnarray*}
So we insert this information into \eqref{TYUAG789Aaasfd}, exchanging the roles
of $x$ and $y$, thanks to the even symmetry of $K$, and we find that
one term symplifies, yielding to the identity
\begin{equation}\label{8iALAJJALA781asA}
\Big({\mathcal{E}}(v;\widetilde{\R}^n)-{\mathcal{F}}(\widetilde v)\Big)
- \Big({\mathcal{E}}(u;\widetilde{\R}^n)-{\mathcal{F}}(u)\Big)
=
\int_{({\widetilde{\R}^n })^c}
\int_{{\widetilde{\R}^n }} 
\widetilde\phi(x)\,\phi(y)
\; K(x - y) \, dx \, dy.\end{equation}
Now, we take a competitor $w$ for $u$ and we define $\phi_1:= (w-u)_+\ge0$
(resp. $\phi_2:= (w-u)_-\ge0$) and $v_1:=u-\phi_1$ (resp. $v_2:=u-\phi_2$). In this way,
we deduce from \eqref{8iALAJJALA781asA} that
\begin{equation} \label{eqaaauatiratysdu}
\Big({\mathcal{E}}(v_i;\widetilde{\R}^n)-{\mathcal{F}}(\widetilde {v_i})\Big)
- \Big({\mathcal{E}}(u;\widetilde{\R}^n)-{\mathcal{F}}(u)\Big)
=
\int_{({\widetilde{\R}^n })^c}
\int_{{\widetilde{\R}^n }}
\widetilde{\phi_i}(x)\,\phi_i(y)
\; K(x - y) \, dx \, dy\ge0,\end{equation}
for $i\in\{1,2\}$.

Notice also that $\max\{ u, w\} = v_1$ and $\min\{ u, w\} = v_2$.
Therefore, after a simple computation we get that
$$ {\mathcal{E}}(v_1;\widetilde{\R}^n)
+{\mathcal{E}}(v_2;\widetilde{\R}^n)\le {\mathcal{E}}(u;\widetilde{\R}^n)
+{\mathcal{E}}(w;\widetilde{\R}^n).$$
This and \eqref{eqaaauatiratysdu} give that
\begin{equation}
\label{u9iocdvbva567a8sd}{\mathcal{E}}(u;\widetilde{\R}^n)
+{\mathcal{E}}(w;\widetilde{\R}^n) \ge
2\,\Big({\mathcal{E}}(u;\widetilde{\R}^n)-{\mathcal{F}}(u)\Big)
+{\mathcal{F}}(\widetilde {v_1})+{\mathcal{F}}(\widetilde {v_2}).
\end{equation}
Now, if $u$ is a 
minimizer of~$\mathcal{F}$ in the class~${\mathcal{A}}_\omega^M$,
we have that ${\mathcal{F}}(u)\le {\mathcal{F}}(\widetilde {v_i})$,
for $i\in\{1,2\}$, and this information, combined with \eqref{u9iocdvbva567a8sd},
gives that ${\mathcal{E}}(w;\widetilde{\R}^n)\ge{\mathcal{E}}(u;\widetilde{\R}^n)$.
This establishes \eqref{9iAJJ1A}.

Now, while the functions in~${\mathcal{M}}_\omega^M$ are 
indeed minimizers of~$\mathcal{E}$ inside~${\widetilde{\mathcal{S}}}_\omega^M$
thanks to \eqref{9iAJJ1A}, there is still no evidence of why they should extend their minimizing properties beyond such domain, and in fact in general they do not. In addition, the set of minimizers~${\mathcal{M}}_\omega^M$ is typically made up of more than just one element. This lack of uniqueness may lead to a corresponding lack of symmetry and rigidity in the elements of~${\mathcal{M}}_\omega^M$, and prevent them from being class~A minimizers.

For these reasons, we direct our attention to a specific element of the set~${\mathcal{M}}_\omega^M$, namely the~\emph{minimal minimizer}. 

The minimal minimizer~$u_\omega^M$ of the class~${\mathcal{M}}_\omega^M$ is defined as
$$
u_\omega^M(x) := \inf_{u \in {\mathcal{M}}_\omega^M} u(x), \qquad \mbox{ for a.e. } x \in \R^n.
$$
It is not hard to prove that~$u_\omega^M$ is unique and belongs to~${\mathcal{M}}_\omega^M$. Therefore, by considering~$u_\omega^M$ we select the element of~${\mathcal{A}}_\omega^M$ with the lowest energy~$\mathcal{F}$ and, at the same time, having an interface with the tightest possible oscillation. As a matter of fact, this double optimality translates into the following two nice features that are enjoyed by the minimal minimizer:
\begin{enumerate}[(i)]
\item the~\emph{doubling} or~\emph{no-symmetry-breaking property}, 
i.e.~$u_\omega^M$ is a minimizer also within functions which exhibit a periodicity of multiple period;
\item the~\emph{Birkhoff property}, i.e.~the level sets of~$u_\omega^M$ and its translations along vectors~$k \in \tau \Z^n$ are well-ordered and have no non-trivial intersections.
\end{enumerate}

The consequences entailed by these facts are twofold. On the one hand, the doubling property implies that~$u_\omega^M$ is a minimizer for~$\mathcal{E}$ with respect to all compact perturbations occurring inside the strip~${\mathcal{S}}_\omega^M$. On the other hand, the Birkhoff property extends such minimizing character to the whole space~$\R^n$, if the width~$M$ of the strip is sufficiently large.

More precisely, by combining the energy and density estimates 
of Theorems~\ref{EBO} and~\ref{DENS}, we obtain the existence 
of a~``clean'' ball~$B$ inside~${\mathcal{S}}_\omega^M$ over 
which~$u_\omega^M$ is, say, smaller than~$- 9/10$. 
The scale invariance of such estimates ensures that the radius 
of this ball is a universal fraction of~$M$. Then, 
the Birkhoff property allows us to translate~$B$ around~${\mathcal{S}}_\omega^M$ 
(in a discrete way) and clean out a full substrip of width comparable to~$\tau$, 
provided that~$M \ge M_0 \tau$, for some large universal constant~$M_0 > 0$. 
This says that the minimal minimizer~$u_\omega^M$ starts attaining values below~$- 9/10$ well before meeting the upper constraint~$\{ \omega \cdot x = M |\omega| \}$. Such~\emph{unconstrainedness} is the key observation that leads to deducing that~$u_\omega^M$ is a class~A minimizer for~$\mathcal{E}$, which thereafter follows almost immediately.

The proof of Theorem~\ref{PLthm} is essentially complete for 
the case of a direction~$\omega \in \tau \Q^n \setminus \{ 0 \}$. 
When instead~$\omega \in \R^n \setminus \Q^n$, we consider a 
sequence~$\{ u_k \}$ of periodic planelike minimizers, 
corresponding to rational approximating directions~$\{ \omega_k \} \subset \tau \Q^n \setminus \{ 0 \}$ of~$\omega$. Uniform H\"older estimates combined with the fact that the value~$M_0$ does not depend on the chosen direction then allow us to take the limit in~$k$ and obtain a planelike class~A minimizer with interface confined in a strip orthogonal to~$\omega$.
\medskip

In conclusion, with Theorem~\ref{PLthm} we are able to prove the existence of planelike minimizers for the energy~$\mathcal{E}$, even if its nonlocal nature prevents it from matching the underlying periodic medium (recall~\eqref{EneE}) as well as in the classical, local case treated in~\cite{V-CRELLE}.
\medskip

On top of that, the techniques described in this section 
are flexible enough to be adapted to other nonlocal periodic models, 
which arise, for instance, in connection with fractional perimeter functionals and long-range Ising models (see~\cite{cozzi-valdinociB, ising}).

\section{Multiwell potentials and chaotic orbits}

Among the many others, one of the main achievements of Paul Rabinowitz
consists in providing a clear and elegant framework
in which the chaotic behavior of 
many equations of great physical importance 
can be rigorously detected and deeply understood.
We would like to point out one application of the
theory that he developed in the framework of Hamiltonian 
dynamics to nonlocal equations. For this,
we consider an equation of the tipe
\begin{equation}\label{D1-MULT}
(-\Delta)^s u(x)+a(x)\,V'(u(x))=0, \qquad{\mbox{ for any }}x\in\R.
\end{equation}
Here, $V$ is a smooth multiwell potential
with a discrete set of nondegenerate minima and~$a$ is
a smooth function (indeed, 
a more general setting, also comprising systems of equations,
may be taken into account).
The reader may compare equations~\eqref{Ds}
and~\eqref{D1-MULT}: in some sense,
equation~\eqref{D1-MULT} is taking into account the possibility
of drifting from one minimum of~$V$ to another one
(as well as the layer solution~$u_0$
in~\eqref{trale} connects the minima of
the two-well potential~$W$).
In this sense, the modulation function~$a$
provides the possibility of favoring this kind of multiple jumps.
\medskip

Equation~\eqref{D1-MULT} also arises naturally in the study
of crystal dislocation dynamics: in this framework,
the function~$u$ can be interpreted as the discrepancy between
the rest position of an atom and its actual
position, see e.g.~Section~2 of~\cite{phy}
for simple physical motivations.\medskip

The result that we present here has been obtained in~\cite{DPV-Rabi}
and can be interpreted as a fractional counterpart of
the classical results in~\cite{RabCZ}.
\medskip

Roughly speaking, these results aim to provide
a {\em symbolic dynamics} for the solutions of
equation~\eqref{D1-MULT} (under suitable
``nondegeneracy assumptions'' on~$a$).
That is, one considers 
the discrete space consisting of 
the equilibria of~$V$ and finds
a solution of~\eqref{D1-MULT}
which induces
a shift operator on such space. Namely,
given a sequence of minima of~$V$
in a prescribed order, one finds
a solution of~\eqref{D1-MULT}
which gets close to each of these minima of~$V$,
in the required order.
\medskip

More precisely, we have the following result:

\begin{theorem}\label{GH:THa}
Assume that~$V\in C^2(\R,[0,+\infty))$
is even, with~$V(r+k)=V(r)$, for any~$r\in\R$ and~$k\in\Z$.
Assume also that~$V(k)=0$ for any~$k\in\Z$,
that~$V(r)>0$ and~$V''(0)>0$.

Let~$a(x):=a_1+a_2\cos(\eps x)$, with~$a_1>a_2>0$ and~$\eps>0$
small enough.

Let~$\zeta_1,\dots,\zeta_N\in \Z$.
Then, there exist~$b_1<\ldots<b_{2N-2}
\in\R$ and
a solution~$u$ of~\eqref{D1-MULT}
such that
\begin{eqnarray*}
&& \lim_{x\to-\infty} u(x)=\zeta_1,\\
&& \big| u(x)-\zeta_1\big| < \frac{1}{10}\quad  {\mbox{ for any }}
x\in (-\infty,b_1],\\
&& \big| u(x)-\zeta_{i+1}\big| < \frac{1}{10} \quad {\mbox{ for any }}
x\in [b_{2i},b_{2i+1}] {\mbox{ for any }}
i\in\{1,\dots,N-2\},\\
&& \big| u(x)-\zeta_N\big| < \frac{1}{10} \quad {\mbox{ for any }}
x\in [b_{2N-2},+\infty)\\ {\mbox{and }}&&
\lim_{x\to+\infty} u(x)=\zeta_N.
\end{eqnarray*}
\end{theorem}

\begin{figure}
    \centering
    \includegraphics[width=16.8cm]{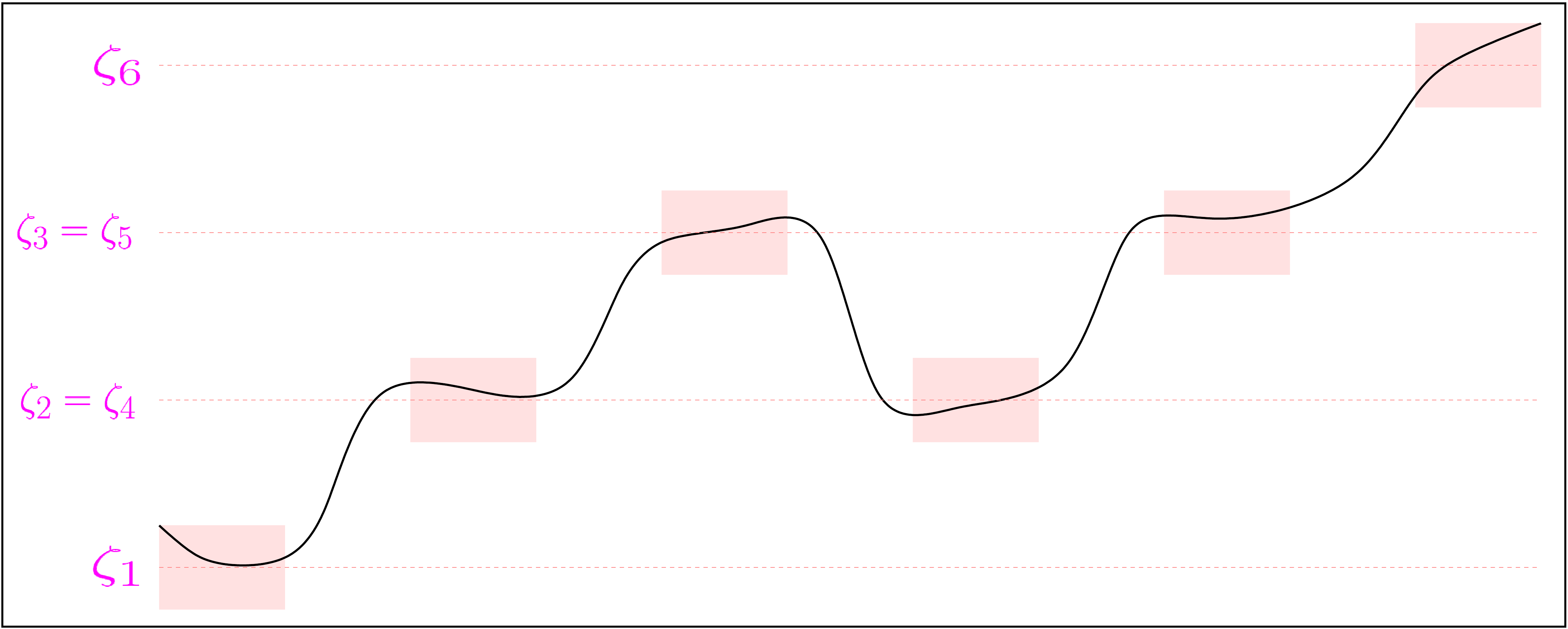}
    \caption{The multibump solution constructed in Theorem~\ref{GH:THa}.}
    \label{fig:bump}
\end{figure}

The situation stated in Theorem~\ref{GH:THa} is described in Figure~\ref{fig:bump}.

Here is a sketch on how
Theorem~\ref{GH:THa} can be deduced from the results in~\cite{DPV-Rabi}.
First of all, we 
may suppose that
\begin{equation}\label{AJ716}
{\mbox{$|\zeta_{i+1}-\zeta_i|=1$
for all~$i\in\{1,\dots,N-1\}$. }}\end{equation}
Indeed, this may be obtained
just by adding the intermediate integers in the original
sequence~$\zeta_1,\dots,\zeta_N$.

Now we claim that
\begin{equation}\label{oi89}
\zeta_{i+1}\in
{\mathcal{A}}(\zeta_i),\end{equation} 
for any~$i\in \{1,\dots,N-1\}$,
where the admissible class~${\mathcal{A}}(\zeta_i)$
is defined in Section~8 of~\cite{DPV-Rabi}
(roughly speaking, ${\mathcal{A}}(\zeta_i)$
contains all the integers~$\zeta$ which can be connected to~$\zeta_i$
by a constrained orbit from a neighborhood of~$\zeta_i$
to a neighborhood of~$\zeta$ with minimal possible action).

To prove~\eqref{oi89}, we suppose that~$\zeta_{i+1}=\zeta_i+1$
(recall~\eqref{AJ716};
the case~$\zeta_{i+1}=\zeta_i-1$ is analogous).
Let~$\zeta \in{\mathcal{A}}(\zeta_i)$.
Since, by Lemma~8.1 in~\cite{DPV-Rabi},
we have that~$2\zeta_i-\zeta$ also belongs to~${\mathcal{A}}(\zeta_i)$,
by possibly replacing~$\zeta$ with~$2\zeta_i-\zeta$
we may and do suppose that~$\zeta>\zeta_i$.
That is, $\zeta\ge\zeta_i+1$. Now, if~$\zeta=\zeta_i+1$, we have that
$$ {\mathcal{A}}(\zeta_i)\ni \zeta=\zeta_i+1=\zeta_{i+1},$$
and so~\eqref{oi89} is proved. 

Hence we may focus on the case in which~$\zeta\ge\zeta_i+2$.
In this case, given~$r\in(0,\,1/4]$,
if we have an orbit~$u_i$ such that~$|u_i(x)-\zeta_i|\le r$
for any~$x\le b_1$ and~$|u_i(x)-\zeta|\le r$
for any~$x\ge b_2$, we can define
$$ u_i^*(x):=\min\{ u_i(x),\,\zeta_i+1\} =\min\{ u_i(x),\,\zeta_{i+1}\}.$$
Then, if~$x\le b_1$, we have that
$$ u_i(x)\le\zeta_i+r < \zeta_i+1,$$
hence~$u_i^*(x)=u_i(x)\in [\zeta_i-r,\,\zeta_i+r]$.
In addition, if~$x\ge b_2$, then
$$ u_i(x)\ge\zeta-r\ge\zeta_i+2-r>\zeta_i+1,$$
thus~$u_i^*(x)=\zeta_i+1 \in [\zeta_{i+1}-r,\,\zeta_{i+1}+r]$.
This gives that~$u^*_i$
is a
constrained orbit from a neighborhood of~$\zeta_i$
to a neighborhood of~$\zeta_{i+1}$.

What is more, $V(u^*_i(x))\le V(u_i(x))$ for any~$x\in\R$;
also, for any~$x$, $y\in\R$,
$$ |u^*_i(x)-u^*_i(y)|\le
|u_i(x)-u_i(y)|.$$
As a consequence, the action of~$u^*_i$ is less than or equal to
the action of~$u_i$, hence~$\zeta_{i+1}$ is admissible,
thus proving~\eqref{oi89}.

Now, from~\eqref{oi89} and Theorem~9.3 in~\cite{DPV-Rabi}
we obtain Theorem~\ref{GH:THa} here.
\medskip

We underline that an important difference between
the result in Theorem~\ref{GH:THa}
and the classical ones for
Hamiltonian systems lies in the glueing methods.
Indeed, in the classical case, the technique of
cutting-and-pasting different trajectories is abundantly used
(tipically, to construct suitable competitors for lowering
the action functional). In the nonlocal case,
this method may lead to additional difficulties,
since the action of the new orbit obtained by a cut-and-paste
of two trajectories
is not simply the sum of the
two actions of the original trajectories, since
nonlocal interactions take place in the elastic
part of the functional, which may be in fact the dominant
contribution when
the fractional parameter~$s$ is small (think once more
to the asymptotics in Theorem~\ref{GAMMA:CO}).
\medskip

To overcome such a difficulty, in~\cite{DPV-Rabi}
we introduced a ``clean interval'' method. Namely,
one has to perform the cut-and-paste techniques
always at points in which the trajectories meet in
a ``very flat'' way (that is the oscillations
of the two trajectories need to be appropriately small
in a sufficiently large interval). This fact,
combined with suitable ``elliptic estimates'',
allows us to estimate the remainder terms and to
efficiently adapt the dynamical systems methods also
to nonlocal cases.

\section*{Acknowledgments}
This work has been supported by 
the ERC grant 277749 ``EPSILON Elliptic PDE's and Symmetry of
Interfaces and Layers for Odd Nonlinearities", the
PRIN grant
201274FYK7 ``Critical Point Theory
and Perturbative Methods for Nonlinear Differential Equations"
and the Alexander von Humboldt Foundation.

\bibliography{AC-bib}

\vspace{2mm}

\end{document}